# ON TWO-DIMENSIONAL AFFINE SPRINGER FIBRES

G. Lusztig

## Introduction

**0.1.** In this paper (until the end of §3) $\Gamma$ is a connected graph with a finite set $I$ of vertices such that for any $i \neq i'$ in $I$ the number $i : i'$ of edges joining $i, i'$ is 0 or 1. (When $i = j$, $i : j$ is not defined.) For $i \in I$ let $m(i) = |\{j \in I; i : j = 1\}|$. We say that $i \in I$ is extremal for $\Gamma$ if $m(i) = 1$. We say that $i$ is a branch point of $\Gamma$ if $m(i) \geq 3$.

**0.2.** We now assume that $\Gamma$ is either

(a) the Coxeter graph of a simple Lie algebra $\mathfrak{g}$ over $\mathbf{C}$ of type $D_n, n \geq 6$ or $E_7, E_8$, or

(b) the affine Coxeter graph of $\mathfrak{g}(K)$ ($K = \mathbf{C}(\epsilon)$ with $\epsilon$ an indeterminate), where $\mathfrak{g}$ is as in (a).

Let $\mathcal{B}$ be the variety of Borel subalgebras of $\mathfrak{g}$ (in case (a)) or the variety of Iwahori subalgebras of $\mathfrak{g}(K)$ (in case (b)).

Let $\mathcal{N}$ be the set of nilpotent elements of $\mathfrak{g}$ (in case (a)) or the set of regular semisimple topologically nilpotent (in the sense of [KL]) elements of $\mathfrak{g}(K)$ (in case (b)).

For $x \in \mathcal{N}$ let $\mathcal{B}_x$ be the set of all $\mathfrak{b} \in \mathcal{B}$ such that $x \in \mathfrak{b}$. This is naturally a nonempty (finite or countable) union of irreducible projective varieties over $\mathbf{C}$, all of the same dimension $d(x)$. (See [S82] in case (a), [KL] in case (b).) The set of projective varieties in this union is denoted by $\mathcal{K}_x$.

Let $\mathcal{N}_2$ be the set of all $x \in \mathcal{N}$ such that $d(x) = 2$ and such that in case (a), $x$ is distinguished (such $x$ exists for each type in (a)) and in case (b), $x$ is obtained by the induction procedure of [KL,§8] from a distinguished nilpotent element $y \in \mathfrak{g}$ with $d(y) = 2$. (Such $y$ exists for each type in (b) and then $x$ is elliptic.)

**0.3.** In the case 0.2(a), the set $\mathcal{K}_x$ for $x \in \mathcal{N}_2$ was studied by Spaltenstein [S82]. He has discovered a remarkable parametrization of $\mathcal{K}_x$ in terms of subgraphs of type $D$ of $\Gamma$. In this paper we prove (in case 0.2(b)) an affine analogue of Spaltenstein's parametrization of $\mathcal{K}_x$.

**0.4.** For $\Gamma$ as in 0.1, let $W$ be the Coxeter group with generators $s_i, i \in I$ and relations $s_i^2 = 1$ for $i \in I$, $s_i s_j = s_j s_i$ if $i : j = 0$ and $s_i s_j s_i = s_j s_i s_j$ if $i : j = 1$. Let $l : W \to \mathbf{N}$ be the length function.

Let $W_1$ be the set of all $w \in W - \{1\}$ that have a unique reduced expression.

Let $W_2$ be the set of all $w \in W-(W_1 \cup \{1\})$ such that for any reduced expression $w = s_{i_1} s_{i_2} \dots s_{i_k}$

(i) we cannot find three consecutive indices $i_a, i_{a+1}, i_{a+2}$, $(1 \le a \le k-2)$ such that $i_a = i_{a+2}$, $i_a : i_{a+1} = 1$;

(ii) we cannot find three consecutive indices $i_a, i_{a+1}, i_{a+2}$, $(1 \le a \le k-2)$ such that $i_a : i_{a+1} = i_{a+1} : i_{a+2} = i_{a+2} : i_a = 0$;

(iii) we cannot find two consecutive indices $i_a, i_{a+1}$ which are extremal and are both joined with the same branch point.

Let $'W_2$ be the set of all $w \in W$ such that $a(w) = 2$ where $a : W \to \mathbf{N}$ is the usual $a$-function.

The following result is proved in [BF98] for $\Gamma$ as in 0.2(a),(b), excluding type $D_n$, affine type $E_7$ and affine type $D_n, n \ge 6$ (but a similar proof applies for affine type $E_7$).

(a) We have $W_2 = {}'W_2$.

Now [GX23] combined with (a) shows that

(b) $|W_2| = |'W_2| < \infty$ for $\Gamma$ of affine type $E_7$ or $E_8$.

(For $\Gamma$ of affine type $E_8$ the finiteness of $W_2$ was proved earlier in [F97].)

Note that if $\Gamma$ is of affine type $E_6$ then $|W_2| = \infty$. Indeed we can write the elements of $I$ as $1,2,3,4,5,6,7$ with $1,5,7$ extremal, $3$ a branch point and $1:2 = 2:3 = 3:4 = 4:5 = 3:6 = 6:7 = 1$. Then $(s_1 s_3 s_2 s_4 s_3 s_5 s_4 s_6 s_3 s_7 s_6 s_2)^d$, $d = 1,2,3,\dots$ form an infinite subset of $W_2$.

Moreover if $\Gamma$ is of affine type $A_n, n \ge 3$ then $|W_2| = \infty$. For example if $n = 3$ we can write the elements of $I$ as $1,2,3,4$ with $1:2 = 2:3 = 3:4 = 4:1 = 1$. Then $(s_1 s_3 s_2 s_4)^d$, $d = 1,2,3,\dots$ form an infinite subset of $W_2$.

In this paper we describe a connection between $\mathcal{K}_x$ (with $x \in \mathcal{N}_2$ and the set $W_2$. To do so we rely on results of Fan [F97], Bremke-Fan [BF98] and Green-Xu [GX23].

**0.5.** For $J \subset I$ we denote by $\underline{J}$ the full subgraph of $\Gamma$ with set of vertices $J$.

Let $\Gamma_2$ be the set of pairs $(J, J!)$ where $J \subset I$ is such that $\underline{J}$ is of type $D_n, n \ge 2$ and $J!$ is a subset of $J$ with $|J!| = 2$ such that if $|J| \ge 3$ then $J!$ consists of extremal elements of $\underline{J}$ both joined with the branch point of $\underline{J}$; moreover, we require that $J_!$ does not consists of two extremal points of $\Gamma$ both joined with a branch point of $\Gamma$.

(When $|J| \ne 4$, $J!$ is uniquely determined by $J$.)

When $|J| \ge 4$ we denote by $j_J$ the extremal point of $\underline{J}$ which is not contained in $J!$. When $|J| = 3$ we denote by $j_J$ the unique element in $J - J!$.

The set $\Gamma_2$ was introduced by Spaltenstein [S82,10.11] in the case where $\Gamma$ is as in 0.2(a).

**0.6.** In this subsection we assume that $\Gamma$ is acyclic. For $(J, J!) \in \Gamma_2$ we define $\tau_{J,J!} \in W$ as follows. If $|J| = 2$ we have $\tau_{J,J!} = 1$. If $|J| \ge 3$ we have $\tau_{J,J!} = s_{i_k} s_{i_{k-1}} \dots s_{i_1}$ where $\{i_1, i_2, \dots, i_k\} = J - J!$ and $i_k = j_J$, $i_1 : i_2 = i_2 : i_3 = \dots =$ $i_{k-1} : i_k = 1$. Using the acyclicity of $\Gamma$ we see that for any $w \in W_2$ there is a unique $(J, J!) \in \Gamma_2$ and a unique $(J', J'!) \in \Gamma_2$ such that

$$w = \tau_{J,J!} s_{i_0} s_{i'_0} z, w = z' s_{j_0} s_{j'_0} \tau^{-1}_{J',J'!},$$

$$J! = \{i_0, i'_0\}, J'! = \{j_0, j'_0\},$$

$$l(w) = l(\tau_{J,J!} s_{i_0} s_{i'_0}) + l(z) = l(z') + l(s_{j_0} s_{j'_0} \tau^{-1}_{J',J'!})$$

for some $z, z'$ in $W$. Thus we obtain a map

(a) $W_2 \to \Gamma_2 \times \Gamma_2$, $w \mapsto ((J, J!), (J', J'!))$.

From the definition

(b) if $w \in W_2$ is mapped by (a) to $((J, J!), (J', J'!))$ then $w^{-1}$ is in $W_2$ and is mapped by (a) to $((J', J'!), (J, J!))$.

The following refinement of 0.4(b) is proved in [GX23].

**Proposition 0.7.** *If $\Gamma$ is of type $E_7$ or $E_8$, or of affine type $E_7$ or $E_8$, then 0.6(a) with $W_2$ replaced by $'W_2$ is a bijection. (But then 0.6(a) itself is a bijection by 0.4(a).) In particular $|W_2|$ is a square.*

(In the case of affine $E_8$ the last sentence was proved earlier in [F97].)

A proof similar to that in [GX23] gives the following result.

**Proposition 0.8.** *If $\Gamma$ is of type $D_n, n \geq 6$, or affine type $D_n, n \geq 6$, then 0.6(a) is a bijection.*

Note that if $\Gamma$ is of affine type $D_5$, then $|W_2| = \infty$. Indeed, we write the elements of $I$ as $1, 2, 3, 4, 5, 6$ with $1, 2$ extremal joined with 3 and $5, 6$ extremal joined with 4 and $3 : 4 = 1$. Then $(s_1 s_4 s_3 s_5 s_4 s_2 s_3 s_6)^d$, $d = 1, 2, 3, \ldots$ form an infinite subset of $W_2$.

**0.9.** Combining 0.7 or 0.8 with 0.6(b) we see that 0.6(a) induces a bijection

(a) $\{w \in W_2; w^2 = 1\} \xrightarrow{\sim} \Gamma_2$.

For $(J, J!) \in \Gamma_2$ the corresponding involution in $W_2$ is

(b) $w_{J,J!} = \tau_{J,J!} s_{i_0} s_{i'_0} \tau^{-1}_{J,J!}$ where $J! = \{i_0, i'_0\}$.

**0.10.** Let $\Gamma$ be as in 0.2 and let $x \in \mathcal{N}_2$. In 2.9 we define a bijection

(a) $\mathcal{K}_x \xrightarrow{\sim} \Gamma_2$.

(In the case 0.2(a) this recovers a result of [S82].)

In §1 we show that in cases 0.2(a),(b), the **C**-vector space $[\Gamma_2]$ is naturally a $W$-module. In 1.8,1.9,1.10 we show that in cases 0.2(b), the restriction of the $W$-module $[\Gamma_2]$ to the finite Weyl group **W** admits a canonical filtration compatible with the basis whose terms are **W**-submodules and whose subquotients are irreducible.

In §3 we show that in the cases 0.2(b), the reflection representation of $W$ can be viewed as a $W$-submodule of its second symmetric power and as a consequence $[\Gamma_2]$ admits a filtration compatible with the $W$-action whose subquotients are irreducible as $W$-modules and as $\mathbf{W}$-modules. Moreover this second filtration is opposed to the first filtration.

In §4 we study affine Springer of dimension 1. We also state a conjecture about decomposing an affine Weyl group into subsets indexed by the nilpotent orbits in the corresponding simple Lie algebra. This decomposition is similar to but in general not the same as the decomposition of the affine Weyl group into two-sided cells.

**0.11.** I wish to thank C.K.Fan, R.M.Green, T.Xu and Z.Yun for valuable comments.

## 1. The $W$-module $[\Gamma_2]$

**1.1.** Let $\Gamma, W$ be as in 0.1, 0.4. Let $V$ be the $\mathbf{C}$-vector space with basis $\{e_i; i \in I\}$. It is well known that there is a well defined $W$-module structure on $V$ such that for $i, j$ in $I$ we have $s_i e_j = e_j + [i,j]e_i$ where $[i,j] = i : j$ if $i \neq j$, $[i,i] = -2$ if $i \in I$. (The reflection representation.) Let $(,) : V \times V \to \mathbf{C}$ be the symmetric bilinear form given by $(e_i, e_j) = -[i,j]$ for all $i, j$ in $I$. We have $(wx, wy) = (x, y)$ for all $x, y$ in $V$ and all $w \in V$.

**1.2.** Let $S^2V \subset V \otimes V$ be the second symmetric power of $V$ with the $W$-module structure induced by that of $V$. Now $S^2V$ has a standard basis consisting of $e_i \otimes e_j + e_j \otimes e_i$ for the various unordered pairs $i \neq j$ in $I$ and of $e_i \otimes e_i$ for various $i \in I$. Let $\phi : S^2V \to \mathbf{C}$ be the restriction of the linear function

$V \otimes V \to \mathbf{C}$, $x \otimes y \mapsto (x, y)$;

let $'S^2V = \ker \phi$. This is a codimension 1 subspace of $S^2V$ stable under the $W$-action. We have

$$\dim {}'S^2V = (1/2)|I|(|I| + 1) - 1 = (1/2)(|I| - 1)(|I| + 2).$$

We show that if $i, j$ in $I$ are such that $i : j = 0$ then:

(a) $$\{f \in S^2V; s_i f = -f, s_j f = -f\} = \mathbf{C}(e_i \otimes e_j + e_j \otimes e_i).$$

It is easy to see that

$$\{f \in V \otimes V; s_i f = -f\} = \oplus_{k;k \neq i} \mathbf{C}(e_i \otimes e_k) \oplus \oplus_{k;k \neq i} \mathbf{C}(e_k \otimes e_i).$$

A similar equality holds for $j$ instead of $i$. Hence

$$\{f \in V \otimes V; s_i f = -f, s_j f = -f\} = \mathbf{C}(e_i \otimes e_j) \oplus \mathbf{C}(e_j \otimes e_i)$$

and (a) follows.

**1.3.** One can verify that

(a) if $\Gamma$ is as in 0.2(a),(b) but not of type $D_n, n \geq 6$ or affine type $D_n, n \geq 6$ then $|\Gamma_2| = (1/2)(|I| - 1)(|I| + 2)$ so that $|\Gamma_2| = \dim {}'S^2V$;

(b) if $\Gamma$ is of type $D_n, n \geq 6$ we have $|\Gamma_2| = (1/2)|I|(|I| - 1)$; if $\Gamma$ is of affine type $D_n, n \geq 6$ we have $|\Gamma_2| = (1/2)(|I|^2 - |I| + 2)$.

(The formulas for $|\Gamma_2|$ appears in [S82,10.11] in the cases 0.2(a).)

**1.4.** For $\Gamma$ as in 0.1 let $[\Gamma_2]$ be the $\mathbf{C}$-vector space with basis $\{\mathbf{e}_{J,J!}; (J,J!) \in \Gamma_2\}$. For any $i \in I$ we define a linear map $s_i : [\Gamma_2] \to [\Gamma_2]$ by

$s_i\mathbf{e}_{J,J!} = -\mathbf{e}_{J,J!}$ if $|J| \ge 3$, $i = j_J$, or if $|J| = 2, i \in J$,

$s_i\mathbf{e}_{J,J!} = \mathbf{e}_{J,J!}$ if $i : j = 0$ for all $j \in J$, or if $|J| \ge 3$, $i \in J$, $i : j_J = 0$,

$s_i\mathbf{e}_{J,J!} = \mathbf{e}_{J',J'!} + \mathbf{e}_{J,J!}$

if $i \notin J$, $|J| \ge 3$, $i : j_J = 1$, $J' = J \cup \{i\}$, $J'! = J!$,

or if $i \notin J$, $|J| \ge 3$, $i : j = 1$ for some $j \in J!$, $J' = J'! = \{i, j_J\}$,

or if $i \in J$, $|J| \ge 3$, $i : j_J = 1$, $J' = J - \{j_J\}$, $J'! = J!$,

or if $|J| = 2$, $i : j = 1$ for any $j \in J$, $J' = J \cup \{i\}, J'! = J!$,

or if $|J| = 2$, $i : j = 1$ for exactly one $j \in J$, $J' = J'! = (J \cup \{i\}) - \{j\}$.

(We note that if $(J,J!) \in \Gamma_2$, $|J| \ge 3$, $i \in I - J, j \in J$ and $i : j = 1$, then $j$ is extremal for $J$.)

**1.5.** For $\Gamma$ as in 0.1 we define a linear map $h : [\Gamma_2] \to S^2V$ by

(a) $h(\mathbf{e}_{J,J!}) = e_j \otimes e_{j'} + e_{j'} \otimes e_j$ if $J = J! = \{j, j'\}$,

(b) $h(\mathbf{e}_{J,J!}) = e_{j_J} \otimes \alpha_J + \alpha_J \otimes e_{j_J}$ if $|J| \ge 3$,

where

$$\alpha_J = e_{j_J} + \sum_{j \in J!} e_j + 2 \sum_{j \in J, j \ne j_J, j \notin J!} e_j.$$

From the definition we see that

(c) $h(s_i\mathbf{e}_{J,J!}) = s_i h(\mathbf{e}_{J,J!})$ for any $i \in I, (J,J!) \in \Gamma_2$.

**Proposition 1.6.** *Assume that $\Gamma$ is as in 0.2(a),(b). Then $h$ is injective. In particular, the maps $s_i : [\Gamma_2] \to [\Gamma_2]$ define a $W$-module structure on $[\Gamma_2]$ so that $\Gamma_2]$ is a $W$-submodule of $S^2V$.*

Let $\mu = \sum_{(J,J!) \in \Gamma_2} c_{J,J!} h(\mathbf{e}_{J,J!})$ where $c_{J,J!} \in \mathbf{C}$. Assume that $\mu = 0$.

We first assume that

(a) $\Gamma$ is as in 1.3(a).

(i) If $i \in I$ is terminal, the coefficient of $e_i \otimes e_i$ in $\mu$ comes from a single $(J,J!)$ so that we must have $c_{J,J!} = 0$.

(ii) Now assume that $i$ is not terminal, not the branch point and that $i : j = 1$ for some terminal $j$. Then the coefficient of $e_i \otimes e_i$ in $\mu$ comes from two elements $(J_1, J_1!), (J_2, J_2!)$ of $\Gamma_2$. We have $i : k = 1$ for a unique $k \in I - \{i, j\}$. We can assume that $(J_1, J_1!) = (\{i, k, j\}, \{i, k\})$, $i \in J_2, k \in J_2, j \notin J_2$. Setting $c_1 = c_{J_1,J_1!}, c_2 = c_{J_2,J_2!}$ we have

$$\begin{aligned} 0 &= c_1(e_i \otimes (e_i + e_j + \dots) + (e_i + e_j + \dots) \otimes e_i) \\ &+ c_2(e_i \otimes (e_i + \dots) + (e_i + \dots) \otimes e_i) + \\ &\text{linear combination of } e_a \otimes e_b \text{ with } (a,b) \ne (i,i), (i,j) \end{aligned}$$

where $+\dots$ refers to terms involving only $e_p$ such that $p \notin \{i, j\}$. (We use (i).) Looking at the coefficient of $e_i \otimes e_j$ we see that $c_1 + c_2 = 0, c_1 = 0$. Thus $c_{J_1,J_1!} = c_{J_2,J_2!} = 0$.

(iii) Next we assume that $i$ is not terminal, not the branch point and that it is not as in (ii). Then the coefficient of $e_i \otimes e_i$ in $\mu$ comes from two elements $(J_1, J_1!), (J_2, J_2!)$ of $\Gamma_2$. We have $i : j = 1$ for some $j$ as in (ii) and $i : k = 1$ for some $k$ not as in (ii). We can assume that $(J_1, J_1!) = (\{i, k, j\}, \{i, k\})$, $i \in J_2, k \in J_2, j \notin J_2$.

Setting $c_1 = c_{J_1,J_1!}, c_2 = c_{J_2,J_2!}$ we have

$$0 = c_1(e_i \otimes (e_i + e_j + \dots) + (e_i + e_j + \dots) \otimes e_i)$$
$$+ c_2(e_i \otimes (e_i + \dots) + (e_i + \dots) \otimes e_i)+$$
$$\text{linear combination of } e_a \otimes e_b \text{ with } (a, b) \neq (i, i), (i, j)$$

where $+\dots$ refers to terms involving only $e_p$ such that $p \notin \{i, j\}$. (We use (ii).) Looking at the coefficient of $e_i \otimes e_j$ we see that $c_1 + c_2 = 0, c_1 = 0$. Thus $c_{J_1,J_1!} = c_{J_2,J_2!} = 0$.

Continuing in this way we see that $c_{J,J!} = 0$ whenever $|J| \geq 3$ and $j_J$ is not the branch point.

Now let $i$ be the branch point. Then the coefficient of $e_i \otimes e_i$ in $\mu$ comes from three elements $(J_1, J_1!), (J_2, J_2!), (J_3, J_3!)$ of $\Gamma_2$. We have $J_1 = \{i, j, k\}, J_2 = \{i, k, p\}, J_3 = \{i, j, p\}$ where $j, k, p$ are the three elements of $\Gamma$ that are joined with $i$. Setting $c_1 = c_{J_1,J_1!}, c_2 = c_{J_2,J_2!}, c_3 = c_{J_3,J_3!}$ we have

$$0 = c_1(e_i \to (e_i + e_j + e_k) + (e_i + e_j + e_k) \otimes e_i)$$
$$+ c_2(e_i \to (e_i + e_k + e_p) + (e_i + e_k + e_p) \otimes e_i)$$
$$+ c_3(e_i \to (e_i + e_j + e_p) + (e_i + e_j + e_p) \otimes e_i)+$$
$$\text{a linear combination of } e_a \otimes e_b \text{ with } a : b = 0.$$

It follows that $c_1 + c_2 = c_1 + c_3 = c_2 + c_3 = 0$ hence $c_1 = c_2 = c_3 = 0$.

We now see that $c_{J,J!} = 0$ whenever $|J| \geq 3$. Since the terms $h(\mathbf{e}_{J,J!})$ with $|J| = 2$ are clearly linearly independent we see that $c_{J,J!} = 0$ for any $(J, J!)$. This completes the proof in the case (a).

Next we assume that

(b) $\Gamma$ is as in 1.3(b).

(i) If $i \in I$ is not a branch point, the coefficient of $e_i \otimes e_i$ in $\mu$ comes from a single $(J, J!) \in \Gamma_2$ so that we must have $c_{J,J!} = 0$.

(ii) If $i \in I$ is a branch point, the coefficient of $e_i \otimes e_i$ in $\mu$ comes from two elements $(J_1, J_1!), (J_2, J_2!)$ of $\Gamma_2$. We have $(J_1, J_1!) = (\{i, k, p\}, \{k, p\})$, $(J_2, J_2!) = (\{i, j, p\}, \{j, p\})$, where $j \neq k$ are extremal and $p$ is not a branch point and not extremal (here we use that $n \geq 6$). Setting $c_1 = c_{J_1,J_1!}, c_2 = c_{J_2,J_2!}$ we have

$$0 = c_1(e_i \to (e_i + e_p + e_k) + (e_i + e_p + e_k) \otimes e_i)$$
$$+ c_2(e_i \to (e_i + e_j + e_p) + (e_i + e_j + e_p) \otimes e_i)+$$
$$\text{a linear combination of } e_a \otimes e_b \text{ with } a : b = 0.$$

Here we use (i). Looking at the coefficient of $e_i\otimes e_k$ and $e_i\otimes e_j$ we see that $c_1=c_2=0$. We now see that $c_{J,J!}=0$ whenever $|J|\ge 3$. Since the terms $h''(\mathbf{e}_{J,J!})$ with $|J|=2$ are clearly linearly independent we see that $c_{J,J!}=0$ for any $(J,J!)\in\Gamma_2$. This completes the proof in the case (b).

**1.7.** From the definition we see that the image of $h$ in 1.6 is contained in $'S^2V$. If we are in case 1.3(a) then $h:[\Gamma_2]\to S^2V$ induces an isomorphism

(a) $[\Gamma_2]\xrightarrow{\sim}{'S^2V}$.

(We use that $\dim[\Gamma_2]=\dim{'S^2V}$, see 1.3(a).)

We see also that in case 1.3(a), $h$ defines a bijection of the basis of $[\Gamma_2]$ indexed by $\Gamma_2$ onto a basis of $'S^2V$. As R.M.Green and T.Xu informed us, this basis of $'S^2V$ coincides with the basis of $'S^2V$ defined in [GX25] (although the definition of the basis in [GX25] does not involve the set $\Gamma_2$).

**1.8.** In this subsection we assume that $\Gamma$ is of affine type $E_8$. We write the elements of $I$ as $1,2,\dots,8,9$ where $1:2=2:3=\dots=7:8=1, 3:9=1$ (the branch point is 3; the extremal points are $1,8,9$). Let $\mathbf{W}$ be the (Weyl) subgroup generated by $s_i; i\ne 8$. We want to study the restriction of the $W$-module $[\Gamma_2]$ to $\mathbf{W}$. Let $A$ be the $\mathbf{C}$-subspace of $[\Gamma_2]$ spanned by the basis vectors corresponding to the following $(J,J!)$:

$$(13,13),(14,14),(15,15),(16,16),(17,17),(19,19),$$
$$(24,24),(25,25),(26,26),(27,27),(29,29),(35,35),(36,36),(37,37),$$
$$(46,46),(47,47),(49,49),(57,57),(59,59),(69,69),(79,79),$$
$$(123,13),(234,24),(239,29),(345,35),(456,46),(439,49),(567,57),$$
$$(2934,29),(29345,29),(293456,29),(2934567,29),(4932,49),$$
$$(49321,49),(2349,24).$$

(These are the $(J,J!)$ which do not involve 8. There are 35 of them.) Let $A'$ be the $\mathbf{C}$-subspace of $[\Gamma_2]$ spanned by the basis vectors corresponding to the following $(J,J!)$:

$$(18,18),(28,28),(38,38),(48,48),(58,58),(68,68),(678,68),(98,98).$$

(These are the $(J,J!)$ which do involve 8 but when $|J|\ge 3$ we have $j_J\ne 8$.) There are 8 of them.)

Let $A''$ be the $\mathbf{C}$-subspace of $[\Gamma_2]$ spanned by the basis vector corresponding to $(J,J!)=(29345678,29)$. (This is the unique $(J,J!)$ such that $|J|\ge 3$, $j_J=8$.)

Then $[\Gamma_2]=A\oplus A'\oplus A''$ and $A$, $A\oplus A'$ are $\mathbf{W}$-submodules. Moreover, $A,(A\oplus A')/A,[\Gamma_2]/(A\oplus A')$ are irreducible $\mathbf{W}$-modules.

**1.9.** In this subsection we asssume that $\Gamma$ is of affine type $E_7$. We write the elements of $I$ as $1,2,\dots,7,8$ where $1:2=2:3=\dots=6:7=1, 4:8=1$ (the branch point is 4; the extremal points are $1,7,8$). Let $\mathbf{W}$ be the (Weyl) subgroup generated by $s_i; i \ne 7$. We want to study the restriction of the $W$-module $[\Gamma_2]$ to $\mathbf{W}$. Let $A$ be the $\mathbf{C}$-subspace of $[\Gamma_2]$ spanned by the basis vectors corresponding to the following $(J, J!)$:

$$(13,13),(14,14),(15,15),(16,16),(18,18),$$
$$(24,24),(25,25),(26,26),(28,28),(35,35),(36,36),(38,38),$$
$$(46,46),(58,58),(68,68),(123,13),(234,24),(345,35),$$
$$(348,38),(456,46),(548,58),(3845,38),(38456,38),(5843,58),$$
$$(58432,58),(584321,58),(3458,35).$$

(These are the $(J, J!)$ which do not involve 7; there are 27 of them.) Let $A'$ be the $\mathbf{C}$-subspace of $[\Gamma_2]$ spanned by the basis vectors corresponding to the following $(J, J!)$:

$$(17,17),(27,27),(37,37),(47,46),(57,57),(567,57),(87,87).$$

(These are the $(J, J!)$ which do involve 7 but when $|J| \ge 3$ satisfy $j_J \ne 7$; there are 7 of them.)

Let $A''$ be the $\mathbf{C}$-subspace of $[\Gamma_2]$ spanned by the basis vector corresponding to $(J, J!) = (384567, 38)$. (This is te unique $(J, J!)$ such that $|J| \ge 3$, $j_J = 7$.) Then $[\Gamma_2] = A \oplus A' \oplus A''$ and $A$, $A \oplus A'$ are $\mathbf{W}$-submodules. Moreover, $A, (A \oplus A')/A, [\Gamma_2]/(A \oplus A')$ are irreducible $\mathbf{W}$-modules.

**1.10.** In this subsection we assume that $\Gamma$ is of affine type $D_n, n \ge 6$. We fix an extremal point $j$ of $\Gamma$. Let $\mathbf{W}$ be the subgroup generated by $s_i, i \in I - \{j\}$ (a Weyl group of type $D_n$). Let $A$ be the subspace of $[\Gamma_2]$ spanned by the basis elements corresponding to the various $(J, J!)$ which do not involve $j$ (there are $(1/2)n(n-1)$ such $(J, J!)$). Let $A'$ be the subspace of $[\Gamma_2]$ spanned by the basis elements corresponding to the various $(J, J!)$ which do involve $j$ but satisfy $j_J \ne j$ when $|J| \ge 3$ (there are $n$ such $(J, J!)$). Let $A''$ be the subspace of $[\Gamma_2]$ spanned by the basis element corresponding to the unique $(J, J!)$ such that $|J| \ge 3$ and $j_J = j$. Then $[\Gamma_2] = A \oplus A' \oplus A''$ and $A$, $A \oplus A'$ are $\mathbf{W}$-submodules. Moreover, $A, (A \oplus A')/A, [\Gamma_2]/(A \oplus A')$ are irreducible $\mathbf{W}$-modules.

**1.11.** In this subsection we assume that $\Gamma$ is acyclic and has a unique branch point and exactly three terminal points. We show:

(a) If $M$ is a $W$-submodule of $[\Gamma_2]$ that contains $\mathbf{e}_{J,J!}$ for some $(J, J!) \in \Gamma_2$ with $|J| = 2$, then $M = [\Gamma_2]$.

From our assumption on $\Gamma$ we see that for any $(J', J'!) \in \Gamma_2$ such that $|J'| = 2$ we can find a sequence $(J_t, J_t!)$ in $\Gamma_2$, $t = 1, 2, \dots, k$ with $|J_t| = 2$ for all $t$ such that $(J_1, J_1!) = (J, J!)$, $(J_k, J_k!) = (J', J'!)$ and such that for any $t = 1, 2, \dots, k-1$, we have $J_t = \{i, i'\}, J_{t+1} = \{i, i''\}$ where $i, i', i''$ in $I$ satisfy $i : i' = 0, i : i'' = 0, i'' : i' = 1$. Using this and the formulas in 1.4 we see by induction that $e_{J_t, J_t!} \in M$ for any $t$ and in particular $e_{J',J'!} \in M$. From the formulas in 1.4 it is clear that the elements $e_{J',J'!}$ with $|J'| = 2$ generate the $W$-module $[\Gamma_2]$. This completes the proof of (a).

**1.12.** The results in this subsection are not used in the rest of the paper. Let $\mathcal{H}$ be the Iwahori-Hecke algebra corresponding to $W$, that is the algebra over $\mathcal{A} := \mathbf{Z}[v, v^{-1}]$ ($v$ is an indeterminate) with generators $T_i, i \in I$ and relations $(T_i + v)(T_i - v^{-1}) = 0$ for $i \in I$, $T_i T_j = T_j T_i$ if $i : j = 0$ and $T_i T_j T_i = T_j T_i T_j$ if $i : j = 1$.

Let $[[\Gamma_2]]$ be the free $\mathcal{A}$-module with basis $\{\mathbf{E}_{J,J!}; (J, J!) \in \Gamma_2\}$. For any $i \in I$ we define an $\mathcal{A}$-linear map $T_i : [[\Gamma_2]] \to [[\Gamma_2]]$ by

$T_i \mathbf{E}_{J,J!} = -v^{-1}\mathbf{E}_{J,J!}$ if $|J| \ge 3$, $i = j_J$, or if $|J| = 2, i \in J$,

$T_i \mathbf{E}_{J,J!} = v\mathbf{E}_{J,J!}$ if $i : j = 0$ for all $j \in J$, or if $|J| \ge 3$, $i \in J$, $i : j_J = 0$,

$T_i \mathbf{E}_{J,J!} = \mathbf{E}_{J',J'!} + v\mathbf{E}_{J,J!}$

if $i \notin J$, $|J| \ge 3$, $i : j_J = 1$, $J' = J \cup \{i\}$, $J'! = J!$,

or if $i \notin J$, $|J| \ge 3$, $i : j = 1$ for some $j \in J!$, $J' = J'! = \{i, j_J\}$,

or if $i \in J$, $|J| \ge 3$, $i : j_J = 1$, $J' = J - \{j_J\}$, $J'! = J!$,

or if $|J| = 2$, $i : j = 1$ for any $j \in J$, $J' = J \cup \{i\}, J'! = J!$,

or if $|J| = 2$, $i : j = 1$ for exactly one $j \in J$, $J' = J'! = (J \cup \{i\}) - \{j\}$.

One can show that if $\Gamma$ is as in 0.2(a),(b), the maps $T_i$ define an $\mathcal{H}$-module structure on $[[\Gamma_2]]$. (We omit the proof.) By specializing $v = 1$ we obtain the $W$-module structure 1.7 on $[\Gamma_2]$.

## 2. Study of $\mathcal{B}_x$

**2.1.** We now assume that $\Gamma$ is as in 0.2(a),(b). Let $x \in \mathcal{N}$ and let $\mathcal{K}_x$ be as in 0.2.

We define a function $\mathcal{K}_x \times \mathcal{K}_x \to W$, by

$(X, X') \mapsto \gamma(X, X')$ =relative position of $\beta, \beta'$

where $(\beta, \beta')$ is a generic point of $X \times X'$. (In the case 0.2(a) this definition appears in [St74], [St76].) From the definition we have $\gamma(X, X') = \gamma(X', X)^{-1}$. In particular,

(a) $\gamma(X, X)$ is an involution.

**2.2.** If $i \in I$, an $i$-line in $\mathcal{B}$ is the set of all $\beta \in \mathcal{B}$ contained in a fixed parabolic (or parahoric) subalgebra of type $i$.

We show:

(a) *In the setup of 2.1, let $i, j$ in $I$ be such that $i : j = 1$. Assume that some $j$-line contained in $\mathcal{B}_x$ meets two different $i$-lines contained in $\mathcal{B}_x$. Then $d(x) \ge 3$.*
Let $\mathfrak{p}$ be the parabolic (in case 0.2(a)) or parahoric (in case 0.2(b)) subalgebra of $\mathfrak{g}$ or $\mathfrak{g}(K)$ generated by the various $\mathfrak{b}$ on the three lines in (a). The reductive quotient $\bar{\mathfrak{p}}$ of $\mathfrak{p}$ is of type $A_2$ and the image $\bar{x}$ of $x$ in $\bar{\mathfrak{p}}$ is a nilpotent element contained in three different maximal parabolic subalgebras. It follows that any $\mathfrak{b}' \in \mathcal{B}$ that is contained in $\mathfrak{p}$ contains $x$ so that $\mathcal{B}_x$ contains a 3-dimensional subvariety. This proves (a).

**2.3.** Let $x \in \mathcal{N}_2$. Let $X \in \mathcal{K}_x, X' \in \mathcal{K}_x$ and let $w = \gamma(X, X') \in W$. Let $w = s_{i_1} s_{i_2} \dots s_{i_k}$ be a reduced expression. We show that

(a) conditions 0.4(i), 0.4(ii) hold.

We can find $\mathfrak{b} \in X, \mathfrak{b}' \in X'$ such that $\mathfrak{b}, \mathfrak{b}'$ are in relative position $w$. There is a unique sequence $\mathfrak{b}_0, \mathfrak{b}_1, \dots, \mathfrak{b}_k$ in $\mathcal{B}$ such that

(b) $\mathfrak{b}_0 = \mathfrak{b}, \mathfrak{b}_k = \mathfrak{b}'$ and $\mathfrak{b}_r, \mathfrak{b}_{r+1}$ are in relative position $s_{i_{r+1}}$ for $r = 0, 1, \dots, k-1$.

We have necessarily $\mathfrak{b}_r \in \mathcal{B}_x$ for $r = 0, 1, \dots, k$. Assume now that we can find $a$ as in 0.4(i) or (ii). Let $\mathfrak{p}$ be the parabolic (in case 0.2(a)) or parahoric (in case 0.2(b)) subalgebra of $\mathfrak{g}$ or $\mathfrak{g}(K)$ generated by $\mathfrak{b}_{i_a}, \mathfrak{b}_{i_{a+1}}$. The reductive quotient $\bar{\mathfrak{p}}$ of $\mathfrak{p}$ is of type $A_2$ or $A_1 \times A_1 \times A_1$; we have $x \in \mathfrak{p}$ and the image $\bar{x}$ of $x$ in $\bar{\mathfrak{p}}$ is a nilpotent element contained in two opposed Borel subalgebras of $\bar{\mathfrak{p}}$ hence $\bar{x} = 0$. It follows that any $\mathfrak{b}' \in \mathcal{B}$ that is contained in $\mathfrak{p}$ contains $x$ so that $\mathcal{B}_x$ contains a 3-dimensional subvariety. This contradicts the assumption that $d(x) = 2$. This proves (a).

We show:

(c) condition 0.4(iii) holds.

Assume the contrary. Then $\mathfrak{b}_{i_a}, \mathfrak{b}_{i_{a+1}}$ are in relative position $s_{i_a} s_{i_{a+1}}$. Let $\mathfrak{p}$ be the parabolic (in case 0.2(a)) or parahoric (in case 0.2(b)) subalgebra of $\mathfrak{g}$ or $\mathfrak{g}(K)$ generated by $\mathfrak{b}_{i_a}, \mathfrak{b}_{i_{a+1}}$. The reductive quotient $\bar{\mathfrak{p}}$ of $\mathfrak{p}$ is of type $A_2$ or $A_1 \times A_1 \times A_1$; we have $x \in \mathfrak{p}$ and the image $\bar{x}$ of $x$ in $\bar{\mathfrak{p}}$ is 0. It follows that $x$ is obtained by induction from the element 0 in $\bar{\mathfrak{p}}$ (in case 0.2(b) induction is in the sense of [KL,§8]). But then $x$ is not distinguished (in case 0.2(a)) and not elliptic (in case 0.2(b)). This contradicts $x \in \mathcal{N}_2$ and proves (c).

**2.4.** We preserve the setup of 2.3. We have $k \geq 1$. (If $k = 0$ then $X = X'$ is a point which is not the case.)

We can find an open dense subset $U'$ of $X'$ and an open dense subset $U$ of $X$ such that for any $\mathfrak{b} \in U, \mathfrak{b}' \in U'$ the relative position of $\mathfrak{b}, \mathfrak{b}'$ is $w$. Let $\mathfrak{b}' \in U'$. For any $\mathfrak{b} \in U$ there is a unique sequence $\mathfrak{b} = \mathfrak{b}_0, \mathfrak{b}_1, \dots, \mathfrak{b}_k = \mathfrak{b}'$ in $\mathcal{B}$ such that 2.3(b) holds. We have necessarily $\mathfrak{b}_r \in \mathcal{B}_x$ for $r = 0, 1, \dots, k$.

By 2.3(a) we can assume that either $i_1 : i_2 \neq 0$ or that $i_1 : i_2 = 0, i_2 : i_3 = 1$.

If $\mathfrak{b}_* \in \mathcal{B}$ belongs to the $i_1$-line through $\mathfrak{b}$ and is different from $\mathfrak{b}$ but is very close to $\mathfrak{b}$ then $\mathfrak{b}_* \in \mathcal{B}_x$ (since the entire $i_1$-line is contained in $\mathcal{B}_x$) and $\mathfrak{b}_* \in U$. Since this holds for any $\mathfrak{b} \in U$ we see that $X$ must be a union of $i_1$-lines. We see also that the map $U \to \mathcal{B}_x$, $\mathfrak{b} \mapsto \mathfrak{b}_1$ has 1-dimensional fibres (parts of $i_1$-lines) hence its image $U_*$ has dimension 1.

Note that

(a) $k \geq 2$.

(If $k = 1$ then $X = X'$ is an $i_1$-line so $d(x) = 1$ which is not the case.)

We are in one of the following cases.

(i) $k = 2$ and $i_1 : i_2 = 0$ or

(ii) $k \geq 2$ and $i_1 : i_2 = 1$ or

(iii) $k \geq 3$ and $i_1 : i_2 = 0$ (in this case using 2.3(a) we can assume after possibly interchanging $i_1, i_2$, that $i_1 : i_3 = 1$).

(In fact we shall eventually show that if $k = 2$ we must be in case (i).) We now assume that we are in case (ii) or (iii). In case (ii) we set $i = i_1, j = i_2$; in case (iii) we set $i = i_1, j = i_3$.

If $\mathfrak{b} \in U$ then the $j$-line through $\mathfrak{b}_1$ is contained in $\mathcal{B}_x$ and any $\mathfrak{b}_{**}$ on this $j$-line from which $\mathfrak{b}_2^{\mathfrak{b}}$ is removed has relative position with $\mathfrak{b}'$ equal to $s_{i_2} s_{i_3} \dots s_{i_k}$. These punctured $j$-lines are disjoint. (If the $j$-line through $\mathfrak{b}_1^{\mathfrak{b}}$ and the $j$-line through $\mathfrak{b}_1^{\tilde{\mathfrak{b}}}$ (with $\tilde{\mathfrak{b}} \in U$ on a different $i$-line) have a common point then $\mathcal{B}_x$ contains a $j$-line which intersects two different $i$-lines. By 2.2(a) this contradicts $d(x) = 2$.) We see that the union of these punctured $j$-lines for various $\mathfrak{b}_{**}$ as above is 2-dimensional (and irreducible). This union must be contained in a unique $X'' \in \mathcal{B}_x$ such that $\gamma(X'', X') = s_{i_2} \dots s_{i_k}$ in case (ii) or $\gamma(X'', X') = s_{i_3} s_{i_4} \dots s_{i_k}$ in case (iii).

We now show that if $k = 2$ then we must be in case (i). (If $i_1 : i_2 = 1$ then by the argument above we can find $X'' \in \mathcal{B}_x$ such that $\gamma(X'', X') = s_{i_2}$, contradicting (a).)

We now study the relationship between $X$ and $X''$ (in case (ii) or (iii)). Let $Z = X \cap X''$. Let $\mathcal{L}$ be the set of $i$-lines contained in $X$. We have an obvious map $Z \to \mathcal{L}$ (any point of $Z$ is contained in a unique $i$-line contained in $X$). From 2.2(a) we see that this map is injective. Thus $Z$ can be identified with a closed subset of the 1-dimensional irreducible projective variety $\mathcal{L}$. Since $Z$ is a 1-dimensional projective variety it follows that $Z \xrightarrow{\sim} \mathcal{L}$. (Similarly any $j$-line in $X''$ contains a unique point in $Z$.) We see that any $i$-line in $\mathcal{L}$ intersects $X''$ and $X$ is the union of all $i$-lines that are contained in $X$ and which intersect $X''$. Now let $\tilde{X} \in \mathcal{K}_x$ be such that $\tilde{X}$ is a union of $i$-lines and such that any $i$-line in $\tilde{X}$ intersects $X''$. We show that $\tilde{X} = X$. Let $L$ be an $i$-line contained in $\tilde{X}$ and let $\mathfrak{b}$ a point in $L \cap X''$. Assume first that $\mathfrak{b} \notin Z$. Let $L'$ be the $j$-line contained in $X''$ such that $\mathfrak{b} \in L'$. As we have seen above, we have $\mathfrak{b}' \in Z$ for some $\mathfrak{b}' \in L'$. Note that $\mathfrak{b}' \neq \mathfrak{b}$ and $\mathfrak{b}'$ is contained in an $i$-line contained in $X$. But this line together with $L, L'$ contradicts 2.2(a). We see that $\mathfrak{b} \in Z$. But the $i$-line through $\mathfrak{b}$ is contained in $X$. Thus $L \subset X$. It follows that $\tilde{X} = X$. We see that

(b) $X \in \mathcal{K}_x$ is characterized by the property that $X$ is a union of $i$-lines and any $i$-line in $X$ intersects $X''$.

**Proposition 2.5.** *Assume that $\Gamma$ is as in 0.2(a),(b). Let $x \in \mathcal{N}_2$. Let $X \in \mathcal{K}_x, X' \in \mathcal{K}_x$ and let $w = \gamma(X, X') \in W$. We have $w \in W_2$.*

By 2.3(a),(c), it is enough to show that $w \notin W_1$. Let $w = s_{i_1} s_{i_2} \dots s_{i_k}$ be a reduced expression. By 2.4 we have $k \geq 2$. We show that $w \notin W_1$ by induction on $k$. If $k = 2$ then from 2.4 we have $i_1 : i_2 = 0$ hence $w \notin W_1$. Assume now that $k \geq 3$. If we are in case (iii) of 2.4 then clearly $w \notin W_1$. If we are in case (ii) of 2.4 then by the argument in 2.4 we can find $X'' \in \mathcal{K}_x$ such that $\gamma(X'', X') = s_{i_2} \dots s_{i_k}$. By the induction hypothesis we have $s_{i_2} \dots s_{i_k} \notin W_1$. It follows that $w \notin W_1$. This completes the proof.

**2.6.** Let $x \in \mathcal{N}_2, X \in \mathcal{K}_x$. Then $\gamma(X, X)$ is an involution in $W_2$ (see 2.1(a),2.5). By 0.9 this equals $w_{J,J!}$ for a well defined $(J, J!) \in \Gamma_2$ (see 0.9(b)); this defines a map $\psi : \mathcal{K}_x \to \Gamma_2$, $X \mapsto (J, J!)$.

**Proposition 2.7.** *Let $x \in \mathcal{N}_2$. Let $X, X'$ in $\mathcal{K}_x$ be such that $\psi(X) = \psi(X') = (J, J!)$. We have $X = X'$. In other words, the map $\psi : \mathcal{K}_x \to \Gamma_2$ is injective.*

We argue by induction on $|J|$. Assume first that $|J| = 2$, $J = \{i, j\}$. Since the relative position of almost any two points in $X$ is $s_i s_j$ we see that $X$ is a union of $i$-lines and a union of $j$-lines. Let $w = \gamma(X, X') \in W_2$. We show that $l(s_i w) = l(w) - 1$. Assume that $l(s_i w) = l(w) + 1$. Now $w$ is the relative position of a generic point $\beta$ of $X$ with a generic point $\beta'$ of $X'$. If $\beta_1 \in X$ is on the same $i$-line as $\beta$ and is very close to $\beta$ (but not equal to $\beta$) then $\beta_1$ is still a generic point of $X$ and the relative position of $\beta_1, \beta'$ is $s_1 w$ (hence not $w$). This contradiction shows that $l(s_i w) = l(w) - 1$. Similarly we have $l(s_j w) = l(w) - 1$. We see that $w = s_i s_j z$ for some $z$ such that $l(w) = l(s_i s_j) + l(z)$. Similarly we have $w = z' s_i s_j$ for some $z'$ such that $l(w) = l(z') + l(s_i s_j) + l(z)$. We see that under the bijection 0.7 or 0.8, $w$ is mapped to $((J_1, J_1!), (J_1, J_1!))$. But $s_i s_j \in W_2$ is mapped by 0.7 or 0.8 to the same element. Using 0.7 or 0.8, we see that $w = s_i s_j$. Using the definition of $w$ we now see that $X = X'$. This completes the proof in our case.

Assume now that $|J| \ge 3$. Let $J_1 = J - \{j_J\}$. We have $|J_1| = |J| - 1$. By 2.4 there exists $X_1 \in \mathcal{K}_x$ such that $\psi(X'') = (J_1, J_1!)$ and such that $X$ is the union of all $j_J$-lines which meet $X_1$. Similarly there exists $X_1' \in \mathcal{K}_x$ such that $\psi(X_1') = (J_1, J_1!)$ and such that $X'$ is the union of all $j_J$-lines which meet $X_1'$. By the induction hypothesis we have $X_1 = X_1'$. It follows that $X = X'$. The proposition is proved.

**Proposition 2.8.** *Let $x \in \mathcal{N}_2$. There exists $X \in \mathcal{K}_x$ such that $\psi(X) = (J, J!)$ where $|J| = 2$.*

Since $\mathcal{K}_x$ is nonempty we can find $X \in \mathcal{K}_x$. We have $\psi(X) = (J, J!)$ for some $(J, J!) \in \Gamma_2$. We argue by induction on $|J|$. If $|J| = 2$ we are done. Assume now that $|J| \ge 3$. Let $J' = J - \{j_J\}$. As in the proof of 2.7 we can find $X' \in \mathcal{K}_x$ such that $\psi(X') = (J', J'!)$ for some $(J', J'!) \in \Gamma_2$. Using the induction hypothesis we see that the proposition holds.

We state the following result which for $\Gamma$ as in 0.2(a) is due to Spaltenstein [S82].

**Theorem 2.9.** *Let $x \in \mathcal{N}_2$. The map $\psi : \mathcal{K}_x \to \Gamma_2$ is a bijection.*

**2.10.** We now prove 2.9 assuming that $\Gamma$ is as in 1.3(a). The set $\mathcal{K}_x$ can be regarded as basis of the homology space $H_4(\mathcal{B}_x, \mathbf{C})$ on which $W$ acts by Springer's representation (in case 0.2(a)) or by the affine version [L96] of the Springer's representation (in case 0.2(b)). Now $W$ also acts naturally on $H_4(\mathcal{B}, \mathbf{C})$ which may be identified with $S^2 V$ with its $W$-module structure in 1.2. The inclusion $\mathcal{B}_x \subset \mathcal{B}$ induces a linear map

(a) $H_4(\mathcal{B}_x, \mathbf{C}) \to H_4(\mathcal{B}, \mathbf{C}) = S^2 V$

compatible with the $W$-module structures. Let $M' \subset S^2 V$ be the image of this map. By 2.8 we can find $X \in \mathcal{K}_x$ such that $\psi(X) = (J, J!)$ where $J = \{i, j\}$, $i : j = 0$. As in the proof of 2.7, $X$ is a union of $i$-lines and a union of $j$-lines. It follows that the image of $X$ in $H_4(\mathcal{B}, \mathbf{C})$ is nonzero and that in $H_4(\mathcal{B}_x, \mathbf{C})$ we have $s_i X = -X, s_j X = -X$. Hence the image $\tau$ of $X$ in $S^2V$ is nonzero and $s_i\tau = -\tau, s_j\tau = -\tau$. By 1.2(a) we have $\tau = c(e_i \otimes e_j + e_j \otimes e_i)$ for some $c \in \mathbf{C}^*$. We have also $\tau \in {}'S^2V$. We see that

$e_i \otimes e_j + e_j \otimes e_i \in M' \cap {}'S^2V$.

Hence the $W$-submodule $M$ of $[\Gamma_2]$ corresponding to $M' \cap {}'S^2V$ contains $\mathbf{e}_{J,J!}$ for some $(J, J!) \in \Gamma_2$ with $|J| = 2$. Using 1.11(a) we see that $M = [\Gamma_2]$ so that $M' \cap {}'S^2V = {}'S^2V$ and $\dim(M' \cap {}'S^2V) = \dim {}'S^2V = |\Gamma_2|$. Now $M'$ is the image of $H_4(\mathcal{B}_x, \mathbf{C})$ (of dimension $|\mathcal{K}_x|$) under a linear map. It follows that

(b) $$|\mathcal{K}_x| \geq \dim(M') \geq \dim(M' \cap {}'S^2V) = \dim {}'S^2V = |\Gamma_2| \geq |\mathcal{K}_x|$$

where the last inequality follows from 2.7. We see that all inequalities in (b) are equalities. Thus we have $M' = {}'S^2V$, $\dim(M') = |\mathcal{K}_x|$ so that (a) restricts to an isomorphism $H_4(\mathcal{B}_x, \mathbf{C}) \xrightarrow{\sim} {}'S^2V$ which can be viewed as an isomorphism

(c) $H_4(\mathcal{B}_x, \mathbf{C}) \xrightarrow{\sim} [\Gamma_2]$.

(It can be shown that this isomorphism maps the basis $\mathcal{K}_x$ of $H_4(\mathcal{B}_x, \mathbf{C})$ onto the basis $\Gamma_2$ of $[\Gamma_2]$; we omit the proof.)

Since $|\mathcal{K}_x| = |\Gamma_2|$ we see from 2.7 that 2.9 holds.

**2.11.** We now prove 2.9 in the case where $\Gamma$ is of affine type $D_n, n \geq 6$. At the same time we give an alternative proof of 2.9 with $\Gamma$ of affine type $E_7$ or $E_8$. In these cases $x \in \mathcal{N}_2$ comes from an element $y \in \mathfrak{g}$ as in 0.2.

It follows from [Y26] that

(a) $$|\mathcal{K}_x| \geq \dim H_4(\mathcal{B}'_y, \mathbf{C}) + \dim H_2(\mathcal{B}'_y, \mathbf{C}) + \dim H_0(\mathcal{B}'_y, \mathbf{C})$$

where $\mathcal{B}'_y$ is the variety of Borel subalgebras of $\mathfrak{g}$ that contain $y$. Hence we have

$|\mathcal{K}_x| \geq 35 + 8 + 1 = 44$ if $\mathfrak{g}$ is of type $E_8$,

$|\mathcal{K}_x| \geq 27 + 7 + 1 = 35$ if $\mathfrak{g}$ is of type $E_7$,

$|\mathcal{K}_x| \geq (1/2)n(n-1) + n + 1 = (1/2)(n^2 + n + 2)$ if $\mathfrak{g}$ is of type $D_n, n \geq 6$,

On the other hand we have

$|\Gamma_2| = 44$ if $\mathfrak{g}$ is of type $E_8$,

$|\Gamma)2| = 35$ if $\mathfrak{g}$ is of type $E_7$,

$|\Gamma_2| = (1/2)(n^2 + n + 2)$ if $\mathfrak{g}$ is of type $D_n, n \geq 6$,

so that $|\mathcal{K}_x| \geq |\Gamma_2|$. Using now 2.7 we see that Theorem 2.9 holds in our case.

We also see that (a) is an equality in our case, confirming a conjecture in [L17, last paragraph of comments to [80]].

**2.12.** We now prove 2.9 in the case where $\Gamma$ is of type $D_n, n \geq 6$. Since it is known that $|\mathcal{K}_x| = |\Gamma_2| = (1/2)n(n-1)$, the result follows from 2.7. A similar argument gives an alternative proof of 2.9 for $\Gamma$ of type $E_7$ or $E_8$.

**Corollary 2.13.** *Let $x \in \mathcal{N}_2$. The map $\gamma : \mathcal{K}_x \times \mathcal{K}_x \to W_2$, $(X, X') \mapsto \gamma(X, X')$ is a bijection.*

It is enough to prove that the composition of $\gamma$ with the inverse of the bijection $W_2 \to \Gamma_2 \times \Gamma_2$ in 0.7, 0.8 is a bijection. This composition is $\mathcal{K}_x \times \mathcal{K}_x \to \Gamma_2 \times \Gamma_2$ and is a product of two bijections as in 2.9, hence is a bijection.

**Proposition 2.14.** *In the setup of 0.2, let $x \in \mathcal{N}_2$. Let $X \in \mathcal{K}_x$. Let $h \in I$ be such that $X$ is a union of $h$-lines. Then the set of $h$-lines contained in $X$ is isomorphic to $P^1$. Hence $X$ is a $P^1$-bundle over $P^1$. In particular $X$ is nonsingular.*

We can find $X' \in \mathcal{K}_x$ such that $\psi(X') = (J, J!)$ where $J = \{i_0, i'_0\}, i_0 : i'_0 = 0$, see 2.8. Let $w = \gamma(X, X') \in W_2$. We can assume that $w = \tau_{J,J!} s_{i_0} s_{i'_0}$ where $\tau_{J,J!} = s_{i_k} s_{i_{k-1}} \dots s_{i_1}$ is as in 0.6. We argue by induction on $|J| = k + 2$. If $|J| = 2$ then $X = X'$ is a union of $i_0$-lines and a union of $i'_0$-lines and the result is obvious. Assume now that $|J| \geq 3$. We associate $X'' \in \mathcal{K}_x$ to $X$ as in 2.4. We have $\psi(X'') = (J'', J!)$ where $J'' = J - \{i_k\}$ and $\gamma(X'', X) = s_{i_{k-1}} \dots s_{i_1} s_{i_0} s_{i'_0}$. Now $X''$ is a union of $i_{k-1}$-lines. Let $Z = X \cap X''$. By the induction hypothesis the set of $i_{k-1}$-lines contained in $X''$ is isomorphic to $P^1$. As shown in 2.4, this set can be identified with $Z$ hence $Z$ is isomorphic to $P^1$. As shown in 2.4, $Z$ can be also identified with the set of $i_k$-lines in $X$. We see that $X$ is a $P^1$-bundle over $P^1$. This completes the proof.

## 3. An imbedding $V \to S^2V$

**3.1.** Let $\mathrm{Hom}_W(V, S^2V)$ be the vector space of linear maps $V \to S^2V$ which commute with the $W$-actions. For a linear map $F : V \to S^2V$ we have $Fe_i = \sum_{jk} a^i_{jk} e_j e_k$ for $i \in I$ where $a^i_{jk} = a^i_{kj} \in \mathbf{C}$ (we omit the $\otimes$ signs). The condition that $F$ commutes with the $W$-actions is

$$F(e_i + [i,p]e_p) = \sum_{jk} a^i_{jk}(e_j + [j,p]e_p)(e_k + [k,p]e_p)$$

for any $i, p$ in $I$ or equivalently

$$\sum_{jk} a^i_{jk} e_j e_k + \sum_{jk} a^p_{jk}[i,p] e_j e_k = \sum_{jk} a^i_{jk} e_j e_k + \sum_{jk} a^i_{jk}[j,p] e_p e_k$$
$$+ \sum_{jk} a^i_{jk}[k,p] e_j e_p + \sum_{jk} a^i_{jk}[j,p][k,p] e_p e_p$$

for any $i, p$ in $I$, or equivalently

$$\sum_{jk} a^p_{jk}[i,p] e_j e_k = \sum_{jk} a^i_{jk}[j,p] e_p e_k$$
$$+ \sum_{jk} a^i_{jk}[k,p] e_j e_p + \sum_{jk} a^i_{jk}[j,p][k,p] e_p e_p$$

for any $i, p$ in $I$, or equivalently

(a) $$a^p_{jk}[i,p] = 0 \text{ for any } p, i, j \neq p, k \neq p,$$

(b) $$[i,p]a^p_{jp} = \sum_k [k,p]a^i_{jk} \text{ for any } p, i, j \neq p,$$

(c) $$a^p_{pp}[i,p] = 2\sum_j a^i_{jp}[j,p] + \sum_{jk} a^i_{jk}[j,p][k,p] \text{ for any } p, i.$$

Now (a) can be rewritten as

(a1) $$a^p_{jk} = 0 \text{ for any } p, j \neq p, k \neq p$$

If $j \neq i$, (b) can be rewritten (using (a1)) as $[i,p]a^p_{jp} = [i,p]a^i_{ji}$. If $i : p = 0$ or $i = p$ this is trivially true so it can be omitted. Thus for $j \neq i$, (b) is equivalent to $a^p_{jp} = a^i_{ji}$ with $i : p = 1, j \neq p, j \neq i$ that is (setting $b^i_j = a^i_{ji} = a^i_{ij}$ when $i \neq j$), to

$$b^p_j = b^i_j \text{ for any } p \neq j, i \neq j, i : p = 1$$

ans also to

(b1) $b^p_j = b^i_j$ for any $j$ and any $i, p$ in the same connected component of $\Gamma - \{j\}$.

If $j = i \neq p$ then (b) can be rewritten (using (a1)) as

$$[i,p]b^p_i = \sum_{k;k\neq i} [k,p]b^i_k + [i,p]a^i_{ii}$$

that is

(b2) $$2b^i_p - \sum_{k;k:p=1} b^i_k = 0 \text{ for any } i, p, i : p = 0$$

(b3) $$b^p_i = -2b^i_p + \sum_{k;k\neq i,k:p=1} b^i_k + a^i_{ii} \text{ for any } i, p, i : p = 1.$$

(In (b2) we have automatically $k \neq i$.) The following statement is a special case of (b2) (assuming (b1)).

(b4) Let $j, l, p$ in $I$ be such that $j : l = 1, l : p = 1, j : p = 0$. We have

$$b^j_l - 2b^l_p + \sum_{k;k:p=1,k\neq l} b^p_k = 0.$$

Indeed $j, l$ are in the same connected component of $\Gamma - \{p\}$ hence $b^j_p = b^l_p$, by (b1); also if $k : p = 1, k \neq l$ then $j, p$ are in the same connected component of $\Gamma - \{k\}$ hence $b^j_k = b^p_k$, by (b1) so that (b4) follows from (b2).

We now show that if (b4) holds then (b2) holds (assuming (b1)). Assume that $i : p = 0$. We can find $j, l$ such that $j, l, p$ are as in (b4) and $i, l$ are in the same connected component of $\Gamma - \{p\}$ (so that $b^l_p = b^i_p$) and $i, p$ are in the same connected component of $\Gamma - \{k\}$ for any $k$ such that $k : p = 1, k \neq l$ (so that $b^i_k = b^p_k$) and $i, j$ are in the same connected component of $\Gamma - \{l\}$ (so that $b^i_l = b^j_l$). We see that (b2) follows from (b4). This (b2) can be replaced by (b4) (assuming (b1)).

If $i = p$, (c) can be rewritten (using (a1)) as

$$-2a^p_{pp} = 2\sum_{j \neq p} b^p_j[j,p] - 4a^p_{pp} + -4\sum_{j \neq p} b^p_j[j,p] + 4a^p_{pp}$$

that is

(c1) $$a^p_{pp} = \sum_{j; j:p=1} b^p_j \text{ for any } p$$

Using (c1) we can rewrite (b3) as

(b5) $$b^p_i + b^i_p = \sum_{k; k \neq i, k:p=1} b^p_k + \sum_{k; k \neq p, k:i=1} b^i_k \text{ for any } i, p, i : p = 1.$$

If $i \neq p$, (c) can be rewritten (using (a1)) as

$$a^p_{pp}[i,p] = 2b^i_p[i,p] + 2\sum_{k \neq i} b^i_k[i,p][k,p] + a^i_{ii}[i,p][i,p].$$

If $i : p = 0$ this is trivially true hence can be omitted. Hence we can assume that $i : p = 1$ and

$$a^p_{pp} - a^i_{ii} = 2b^i_p + 2\sum_{k \neq i, k:p=1} b^i_k - 4b^i_p$$

In the presence of (c1) this is the same as

$$\sum_{j; j:p=1} b^p_j - \sum_{j; j:i=1} b^i_j = -2b^i_p + 2\sum_{j \neq i, j:p=1} b^i_j \text{ with } i : p = 1$$

that is

$$b^p_i + \sum_{j; j \neq i, j:p=1} b^p_j - \sum_{j; j:i=1} b^i_j = -2b^i_p + 2\sum_{j \neq i, j:p=1} b^i_j \text{ with } i : p = 1.$$

Using $b^p_j = b^i_j$, see (b1), this is the same as

$$b^p_i - \sum_{j; j \neq i, j:p=1} b^p_j - \sum_{j; j \neq p, j:i=1} b^i_j - b^i_p = -2b^i_p \text{ with } i : p = 1.$$

This is the same as (b5) hence it can be omitted. We see that the following holds.

**Proposition 3.2.** *Let $E$ be the vector space consisting of all $b_i^p \in \mathbf{C}$ for various $i \neq p$ in $I$ such that (b1),(b4),(b5) hold. Then*

*(a) the map $\mathrm{Hom}_W(V, S^2V) \to E$ described in 3.1 is an isomorphism of vector spaces.*

*(b) For any $F \in \mathrm{Hom}_W(V, S^2V)$ we have $F(V) \subset {}'S^2V$.*

(a) follows from the arguments in 3.1. We prove (b). We must show that $\phi(F(e_i)) = 0$ for any $i \in I$ where $\phi$ is as in 1.2. With notation of 3.1 we have:

$$F(e_i) = a_{ii}^i e_i \otimes e_i + \sum_{j;j\neq i} b_j^i (e_j \otimes e_i + e_i \otimes e_j)$$

hence $\phi(F(e_i)) = 2a_{ii}^i - 2\sum_{j;j:i=1} b_j^i$ and this is 0 by 3.1(c1).

**Corollary 3.3.** *We assume that $\Gamma$ in 0.1 is acyclic. Let $\tilde{E}$ be the vector space consisting of all $b_i^p \in \mathbf{C}$ for various $i \neq p$ in $I$ with $i : p = 1$ such that 3.1(b4),(b5) hold. Then the map $\mathrm{Hom}_W(V, S^2V) \to \tilde{E}$ described in 3.1 is an isomorphism of vector spaces.*

Using the connecteness and acyclicity of $\Gamma$ we can express any $b_j^i, i \neq j$ as a $b_j^p$ where $p : j = 1$ (see 3.1(b1)).) We also use the fact that 3.1(b4),(b5) involve only elements $b_j^i$ with $i : j = 1$.

**3.4.** We now assume that $\Gamma$ is as in 1.11. The elements of $I$ are denoted by $i_1, i_2, \dots, i_p; j_1, j_2, \dots, j_q; k_1, k_2, \dots, k_r$ where $|I| = p + q + r - 2$, $i_1 = j_1 = k_1$ is the branch point,

$i_1 : i_2 = i_2 : i_3 = \dots = i_{p-1} : i_p = 1$, $j_1 : j_2 = j_2 : j_3 = \dots = j_{q-1} : j_q = 1$, $k_1 : k_2 = k_2 : k_3 = \dots = k_{r-1} : k_r = 1$,

and $p \geq 2, q \geq 2, r \geq 2$.

Let $\tilde{E}'$ be the vector space consisting of all $b_i^p \in \mathbf{C}$ for various $i \neq p$ in $I$ with $i : p = 1$ such that 3.1(b4) holds. From the definition we see that $\dim \tilde{E}' = 1$. A nonzero vector in $\tilde{E}'$ can be described as follows.

$$b_{i_{g+1}}^{i_g} = (p-g)/p, g = 1, 2, \dots, p-1,$$

$$b_{i_g}^{i_{g+1}} = (g-1)(1/q + 1/r) - (g-2), g = 1, 2, 3, \dots, p-1,$$

$$b_{j_{g+1}}^{j_g} = (q-g)/q, g = 1, 2, \dots, q-1,$$

$$b_{j_g}^{j_{g+1}} = (g-1)(1/p + 1/r) - (g-2), g = 1, 2, 3, \dots, q-1,$$

$$b_{k_{g+1}}^{k_g} = (r-g)/r, g = 1, 2, \dots, r-1,$$

$$b_{k_g}^{k_{g+1}} = (g-1)(1/p + 1/q) - (g-2), g = 1, 2, 3, \dots, e-1.$$

The condition that this vector satisfies 3.1(b5) is that $1/p + 1/q + 1/r = 1$ or equivalently that $(p, q, r)$ is (up to permutation) $(2, 3, 6), (2, 4, 4)$ or $(3, 3, 3)$, or equivalenty that $\Gamma$ is of affine type $E_8, E_7$ or $E_6$. In these cases our vector satisfies $b_j^i > 0$ for all $i, j$ such that $i : j = 1$.

**3.5.** In this subsection we assume that $\Gamma$ is of affine type $D_n, n \geq 4$. From the definition we see that $\dim \tilde{E} = 1$. A nonzero vector in $\tilde{E}$ can be described as follows. We have $b^i_j = 1$ if $i$ is branch point and $j$ is an extremal point joined with $i$; $b^i_j = 2$ for all other pairs $i, j$ in $I$ with $i : j = 1$.

**3.6.** Assume that $F$ in 3.1 satisfies $F(\sum_i x_i e_i) = 0$ where $x_i \in \mathbf{C}$. We then have

$$\sum_i x_i a^i_{ii} e_i \otimes e_i + \sum_{i,j;j \neq i} x_i b^i_j (e_j \otimes e_i + e_i \otimes e_j) = 0.$$

It follows that $x_i b^i_j = 0$ for all $i \neq j$. If for any $i$ we can find $j$ such that $i : j = 1, b^i_j \neq 0$ then $x_i = 0$ for all $i$ so that $F$ is injective.

**Theorem 3.7.** *(a) If $\Gamma$ is of affine type $E_6, E_7, E_8$ or $D_n, n \geq 4$, then there is a unique (up to a nonzero scalar) linear imbedding $V \to S^2V$ which commutes with the $W$-actions. Its image is contained in ${}'S^2V$.*

*(b) If $\Gamma$ is as in 1.11, but not of affine type $E_6, E_7, E_8$ then any linear map $V \to S^2V$ which commutes with the $W$-actions is* $0$.

This follows from the results in 3.3,3.4,3.5,3.6.

**3.8.** In this subsection we assume that $\Gamma$ in 0.1 is an affine Coxeter graph. We show:

(a) there exists a unique line $L$ in $V$ on which $W$ acts trivially.

Let $u = \sum_{i \in I} a_i e_i \in V$ where $a_i \in \mathbf{C}$ are not all 0. The condition that $s_p u = u$ for all $p \in I$ is that for all $p$ we have $\sum_i a_i e_i = \sum_i a_i (e_i + [i,p] e_p)$ that is $\sum_i a_i [i,p] = 0$. According to Coxeter, for our $\Gamma$, the vector space of functions $\mathcal{C} : I \to \mathbf{C}$ such that $\sum_i [i,k] \mathcal{C}(k) = 0$ for any $k \in I$ is one dimensional. Moreover there is a unique $\mathcal{C}$ in this one dimensional vector space such that $\mathcal{C}(i) \in \mathbf{Z}_{>0}$ for all $i \in I$ and $\mathcal{C}(i) = 1$ for some $i \in I$. We see that $a_i = c\mathcal{C}(i)$ for all $i$ where $c \in \mathbf{C}^*$. This proves (a).

**3.9.** If $\Gamma$ is of affine type $D_n, n \geq 4$ or of affine type $E_6, E_7$, or $E_8$, the 1-dimensional vector space $\tilde{E}$ (see 3.3,3.4,3.5) is spanned by the vector given by $b^i_j = \mathcal{C}(j)$ for all $i, j$ such that $i : j = 1$. In this case the corresponding $a^i_{ii}$ is equal to $2\mathcal{C}(i)$ (see (b3)).

**3.10.** Assume that $\Gamma$ is of affine type $E_7$ or $E_8$. By 1.6 and 3.7(a) we can regard $V$ as a $W$-submodule of $[\Gamma_2]$. Let $\mathbf{W}$ be as in 1.8,1.9. From 1.8,1.9, the $\mathbf{W}$-module $[\Gamma_2]$ is a direct sum of three irreducible $\mathbf{W}$-modules of dimension $35, 8, 1$ (for affine $E_8$) or $27, 7, 1$ (for affine $E_7$). As a $\mathbf{W}$-module, $V$ is a direct sum of two irreducible $\mathbf{W}$-modules of dimension $8, 1$ (for affine $E_8$) or $7, 1$ (for affine $E_7$). Hence the quotient $[\Gamma_2]/V$ viewed as as $\mathbf{W}$-module is irreducible of dimension 35 (for affine $E_8$) or 27 (for affine $E_7$). Note also that $\mathbf{C}$ can be regarded as a $W$-submodule of $V$ by 3.8. We see that $[\Gamma_2]$ has a filtration $0 \subset \mathbf{C} \subset V \subset [\Gamma_2]$ of $[\Gamma_2]$ by $W$-submodules with successive quotients of dimension $1, 8, 35$ (for affine $E_8$) or $1, 7, 27$ (for affine $E_7$). This filtration is opposite to the filtration $0 \subset A \subset A \oplus A' \subset [\Gamma_2]$ by $\mathbf{W}$-modules (compatible with the basis of $[\Gamma_2]$) in 1.8, 1.9, with successive quotients of dimension $35, 8, 1$ (for affine $E_8$) or $27, 7, 1$ (for affine $E_7$).

**3.11.** Assume that $\Gamma$ is of affine type $D_n, n \geq 6$. Let $\mathbf{W}$ be as in 1.10. Now the $\mathbf{W}$-module $'S^2V$ is a direct sum of 4 irreducible modules of dimension $1, n, (1/2)n(n-1), n-1$; its $\mathbf{W}$-submodule $[\Gamma_2]$ is a direct sum of 3 irreducible modules of dimension $1, n, (1/2)n(n-1)$ and its $\mathbf{W}$-submodule $V$ (see 3.7(a)) is a direct sum of 2 irreducible modules of dimension $1, n$. It follows that $V$ is a $W$-submodule of $[\Gamma_2]$ and $[\Gamma_2]/V$ is irreducible of dimension $(1/2)n(n-1)$ for $\mathbf{W}$.

Note also that $\mathbf{C}$ can be regarded as a $W$-submodule of $V$ by 3.8. We see that $[\Gamma_2]$ has a filtration $0 \subset \mathbf{C} \subset V \subset [\Gamma)2]$ by $W$-submodules with successive quotients of dimension $1, n, (1/2)n(n-1)$. This filtration is opposite to the filtration $0 \subset A \subset A \oplus A' \subset [\Gamma_2]$ of $[\Gamma_2]$ by $\mathbf{W}$-modules (compatible with the basis of $[\Gamma_2]$) in 1.10 with successive quotients of dimension $(1/2)n(n-1), n, 1$.

## 4. One dimensional affine Springer fibres

**4.1.** In this section $\Gamma$ denotes the Coxeter graph of $\mathfrak{g}(K)$ where $\mathfrak{g}$ is a simple Lie algebra over $\mathbf{C}$ and $K$ as in 0.2(b).

Then $\mathcal{B}, \mathcal{N}$ are defined as in 0.2. For $x \in \mathcal{N}$, $\mathcal{B}_x, \mathcal{K}_x, d(x)$ are defined as in 0.2. Let $\underline{\mathfrak{g}}$ be a subset of $\mathfrak{g}$ containing exactly one element in each nilpotent orbit in $\mathfrak{g}$. For each $y \in \underline{\mathfrak{g}}$ we denote by $x$ an element of $\mathcal{N} \subset \mathfrak{g}(K)$ obtained from $y$ by the induction procedure [KL,§8]. The elements $x$ so obtained form a subset $\mathcal{N}_*$ of $\mathfrak{g}(K)$ in bijection with $\underline{\mathfrak{g}}$.

Let $W$ be the affine Weyl group of $\mathfrak{g}(K)$.

For $x \in \mathcal{N}_*$ we define $\gamma : \mathcal{K}_x \times \mathcal{K}_x \to W$ as in 2.1. Let $W(x) \subset W$ be its image.

We conjecture that

(a) $$W = \sqcup_{x \in \mathcal{N}_*} W(x).$$

The analogous statement for $\mathfrak{g}$ instead of $\mathfrak{g}(K)$ is known to be true.

**4.2.** In the remainder of this paper we assume that $\Gamma$ is as in 4.1 but not of affine type $A$.

We consider all pairs $(\tilde{\Gamma}, \sigma)$ where $\tilde{\Gamma}$ is a graph as in 4.1 (not of affine type $A$), assumed to be simply laced and $\sigma$ is a graph automorphism of $\tilde{\Gamma}$ of order $\leq 3$ such that

(i) there is no edge or $\tilde{\Gamma}$ that joins two vertices of $\tilde{\Gamma}$ in the same $\sigma$-orbit;

(ii) if $\sigma \neq 1$, the number of extremal points of $\tilde{\Gamma}$ that are fixed by $\sigma$ is 0 or 1.

**4.3.** To any $(\tilde{\Gamma}, \sigma)$ as above we associate a graph $\Gamma$ as in 4.1 as follows. The set of vertices $I$ of $\Gamma$ is the set of orbits of $\sigma$ on the set of vertices $\tilde{I}$ of $\tilde{\Gamma}$. Let $\pi : \tilde{I} \to I$ be the map which associates to any $i \in \tilde{I}$ its $\sigma$-orbit. Then $j \neq j'$ in $I$ are joined in $\Gamma$ whenever $i, i'$ are joined in $\tilde{\Gamma}$ for some $i \in \pi^{-1}(j), i' \in \pi^{-1}(j')$; in that case the number of edges of $\Gamma$ that join $j, j'$ is

$$\max(|\pi^{-1}(j)|/|\pi^{-1}(j')|, |\pi^{-1}(j')|/|\pi^{-1}(j)|.$$

For $i \in \tilde I$ we write $\bar i = \pi(i)$.

Note that $(\tilde\Gamma, \sigma) \mapsto \Gamma$ defines a bijection between the set of pairs $(\tilde\Gamma, \sigma)$ as in 4.2 up to isomorphism and the set of $\Gamma$ as in 4.1 (not of affine type $A$) up to isomorphism.

If $\Gamma$ is of affine type $D_n, n \ge 4$, or $E_6, E_7, E_8$ then $\tilde\Gamma = \Gamma, \tilde I = I, \sigma = 1$.

If $\Gamma$ is of affine type $C_n, n \ge 2$, then $\tilde\Gamma$ is of affine type $D_{n+2}$; we write $\tilde I = \{1, 2, \dots, n+2, n+3\}$; the edges are

$i --- (i+1)$ for $i = 2, 3, \dots, n$ and also $1 --- 3, (n+1) --- (n+2)$, $(n+1) --- (n+3)$;

$\sigma$ acts by

$i \mapsto i$ for $i = 3, \dots, n, n+1$ and $1 \mapsto 2 \mapsto 1$, $n+2 \mapsto n+3 \mapsto n+2$.

We have $I = \{\bar 2, \bar 3, \dots, \overline{n+2}\}$.

If $\Gamma$ is of affine type $B_n, n \ge 3$, then $\tilde\Gamma$ is of affine type $D_{2n}$; we write $\tilde I = \{0, 1, 2, \dots, 2n-1, 2n\}$; the edges are

$i --- (i+1)$ for $i < 2n-1$ and also $0 --- 2, (2n-2) --- 2n$;

$\sigma$ acts by $i \mapsto 2n - i$. We have $I = \{\bar 0, \bar 1, \bar 2, \dots, \bar n\}$.

If $\Gamma$ is of affine type $F_4$ then $\tilde\Gamma$ is of affine type $E_7$; we write

$\tilde I = \{0, 1, 2, 3, 4, 5, 6, 7\}$;

the edges are

$i --- i+1$ for $i = 0, 1, 2, 3, 4$ and $3 --- 6$, $5 --- 7$;

$\sigma$ acts by

$i \mapsto 6 - i$ for $i = 1, 2, 3, 4, 5$, $6 \mapsto 6$, $0 \mapsto 7 \mapsto 0$.

We have $I = \{\bar 0, \bar 1, \bar 2, \bar 3, \bar 6\}$.

If $\Gamma$ is of affine type $G_2$ then $\tilde\Gamma$ is of affine type $E_6$; we write $\tilde I = \{0, 1, 2, 3, 4, 5, 6\}$; the edges are

$1 --- 2$, $2 --- 3$, $2 --- 4$, $0 --- 1$, $4 --- 5$, $3 --- 6$;

$\sigma$ acts by

$2 \mapsto 2, 1 \mapsto 3 \mapsto 4 \mapsto 1, 0 \mapsto 6 \mapsto 5 \mapsto 0$.

We have $I = \{\bar 0, \bar 1, \bar 2\}$.

**4.4.** Let $W_1$ be the set of all $w \in W - \{1\}$ such that $w$ has a unique reduced expression $s_{j_1} s_{j_2} \dots s_{j_k}$ in $W$ and moreover for any $j \ne j'$ in $I$ such that $j : j' \ge 2$ we cannot find $a \in \{1, 2, \dots, k-2\}$ such that $j_a = j, j_{a+1} = j', j_{a+2} = j$ $|\pi^{-1}(j)| < |\pi^{-1}(j')|$. (This is compatible with the definition of $W_1$ in 0.4.)

For example, if $\Gamma$ is of affine type $G_2$ then $W_1$ consists of

$$\bar 0, \bar 1, \bar 2, \bar 0\bar 1, \bar 1\bar 0, \bar 1\bar 2, \bar 2\bar 1, \bar 0\bar 1\bar 2, \bar 2\bar 1\bar 0, \bar 1\bar 2\bar 1, \bar 1\bar 2\bar 1\bar 0, \bar 0\bar 1\bar 2\bar 1, \bar 0\bar 1\bar 2\bar 1\bar 0.$$

(We write $jj' \dots$ instead of $s_j s_{j'} \dots$.)

**4.5.** For $\Gamma$ in 4.1 (not of affine type $A$) we consider $x \in \mathcal N_*$ such that $d(x) = 1$.

Note that the variety $\mathcal B_x$ is a union of $j$-lines in $\mathcal B$ for various $j \in I$. According to [KL,7.7], the set $\mathcal K_x$ is in canonical bijection with $\tilde I$ in such a way that if $i \ne i'$ in $\tilde I$ then the corresponding irreducible components of $\mathcal B_x$ intersect if and only if $i, i'$ are joined in $\tilde{\Gamma}$. One can show that the irreducible component of $\mathcal{B}_x$ corresponding to $i \in \tilde{I}$ is a $j$-line where $j = \pi(i)$.

The following result can be verified in each case.

**Theorem 4.6.** *We have $W(x) = W_1$ (see 4.1). In other words, the map $\mathcal{K}_x \times \mathcal{K}_x \to W$, $(X, X') \mapsto \gamma(X, X')$ (see 4.1) has image equal to $W_1$. Its nonempty fibres are exactly the $\sigma$-orbits on $\mathcal{B}_x \times \mathcal{B}_x$ (identified with $\tilde{I} \times \tilde{I}$).*

Department of Mathematics, M.I.T., Cambridge, MA 02139